\documentclass[11pt]{article}
\usepackage{amsmath}
\allowdisplaybreaks[4]
\usepackage{amsfonts,amsthm,amssymb}
\usepackage{multirow}
\usepackage{array}
\usepackage{makecell}
\usepackage{booktabs}
\usepackage[hidelinks]{hyperref}
\usepackage{graphicx}
\usepackage{graphics}
\usepackage{needspace}

\newtheorem{lem}{Lemma}[section]
\newtheorem{thm}[lem]{Theorem}
\newtheorem{cor}[lem]{Corollary}
\newtheorem{pro}[lem]{Proposition}

\theoremstyle{definition}

\title{The Restricted Edge-Connectivity of Mycielskian Graphs\footnote{The research is supported by Natural Science Foundation of Xinjiang Uygur Autonomous Region (2025D01E02) and National Natural Science Foundation of China (12261086).}}
\author{Huyun Zhu, Yingzhi Tian\thanks{Corresponding author.
E-mail: tianyzhxj@163.com (Y. Tian).}\\
{\small College of Mathematics and System Sciences, Xinjiang University, Urumqi, Xinjiang 830046, PR China}}
\date{}

\begin{document}
\maketitle
\vspace{-1.5em}
\noindent{\bf Abstract} For a graph $G=(V(G),E(G))$, an edge subset
 $F\subseteq E(G)$  is called a
restricted edge-cut if $G-F$ is not connected and  has no isolated vertices. The restricted
edge-connectivity of $G$, denoted by $\lambda'(G)$, is defined as the minimum
cardinality among all restricted edge-cuts if there are at least one; otherwise  
$\lambda'(G)=+\infty$. The Mycielskian graph $\mu(G)$ is defined as the graph
obtained from $G$ by adding shadow vertices $u_i$ corresponding to the vertices
$v_i$ and a root vertex $w$, where each $u_i$ is adjacent to every neighbor of
$v_i$ in $G$ and $w$ is adjacent to every shadow vertex. In this paper, we prove
that $\lambda'(\mu(G))=\min\{\xi(\mu(G)),\beta^{*}(G)\}$, where $\beta^{*}(G)$
is the minimum of $3|\partial_G(X)|+|X|$ over all restricted edge-cuts
$\partial_G(X)$ of $G$ and is set to $+\infty$ if $G$ contains no restricted edge-cut. Moreover,
we show that if $\lambda'(G)=\xi(G)$, then $\lambda'(\mu(G))=\xi(\mu(G))$, where $\xi(G)$ and $\xi(\mu(G)$ are the minimum edge-degrees of $G$ and $\mu(G)$, respectively.


\section{Introduction}\label{sec:intro}

An interconnection network can be modeled by a graph $G=(V(G),E(G))$, where vertices represent processors and edges represent communication links. An edge subset $F\subseteq E(G)$ is called an edge-cut if $G-F$ is disconnected. The edge-connectivity  $\lambda(G)$ of $G$ is the minimum cardinality among all edge-cuts. It is well known that $\lambda(G)\le\delta(G)$, where $\delta(G)$ denotes the minimum degree of $G$. So  $G$ is called maximally edge-connected if $\lambda(G)=\delta(G)$. Furthermore, if every minimum edge-cut of $G$ isolates a vertex,   $G$ is called super edge-connected.

Edge-connectivity imposes no condition on the components that remain after edge deletion. In particular, a minimum edge-cut may isolate a single vertex. To take the structure of the remaining components into account, Esfahanian and Hakimi \cite{esfahanian} introduced the concept of restricted edge-connectivity, which is a kind of conditional edge-connectivity proposed by Harary \cite{harary}. An edge-cut $F$ of $G$ is called a restricted edge-cut if every component of $G-F$ has at least two vertices, equivalently, if $G-F$ has no isolated vertices. The restricted edge-connectivity of $G$, denoted by $\lambda'(G)$, is the minimum cardinality among all restricted edge-cuts if there are at least one; otherwise  
$\lambda'(G)=+\infty$. The degree of a vertex $u$ in $G$ is denoted by $d_G(u)$. For an edge $xy\in E(G)$, its edge-degree is $d_G(x)+d_G(y)-2$. The minimum edge-degree of $G$ is $\xi(G)=\min\{d_G(x)+d_G(y)-2:xy\in E(G)\}$. Esfahanian and Hakimi \cite{esfahanian} proved that every connected graph of order at least four other than a star admits a restricted edge-cut and satisfies $\lambda'(G)\le\xi(G)$. Accordingly, a graph $G$ is called maximally restricted edge-connected if $\lambda'(G)=\xi(G)$.

A useful structural tool for the restricted edge-connectivity \(\lambda'(G)\) is due to Xu and Xu \cite{xu2002}, who proved that if a connected graph $G$ is not maximally restricted edge-connected, then any two distinct \(\lambda'\)-atoms of $G$ are disjoint. Ueffing and Volkmann \cite{ueffing2003} showed that a graph $G$ on $n$ vertices satisfying \(\delta(G)\ge\lfloor n/2\rfloor+1\) is maximally restricted edge-connected. Subsequently, 
Wang and Li \cite{wangli} proved   that $G$ is maximally restricted edge-connected  by relaxing the minimum degree condition above to the Ore-type  condition:   $d_G(u)+d_G(v)\ge n+1$ for every pair of non-adjacent vertices $u$ and $v$.

The restricted edge-connectivity of standard graph products has also been investigated. For Cartesian product, Wang, Ou and Zhu \cite{wangou2010} gave an exact expression for the restricted edge-connectivity of the Cartesian product of two regular graphs; Yokota \cite{yokota} proved that if $G_1$ and $G_2$ are both maximally edge-connected and maximally restricted edge-connected and satisfy $\delta(G_i)\ge2$ for $i=1,2$, then their Cartesian product is maximally restricted edge-connected. For direct product, Ma, Wang and Zhang \cite{ma} obtained an exact expression for the restricted edge-connectivity of the direct product of a nontrivial connected graph $G$ and a complete graph $K_n$. For strong product, Ye and Tian \cite{yetian} determined the restricted edge-connectivity of the strong product of a nontrivial connected graph $G$ with a path $P_n$, a cycle $C_n$, or a complete graph $K_n$.

Another classical graph operation is the Mycielski construction, introduced by Mycielski \cite{mycielski} in 1955 to construct triangle-free graphs with arbitrarily large chromatic number. The Mycielskian graph $\mu(G)$ is obtained from $G$ by adding a shadow vertex corresponding to each vertex of $G$, joining each shadow vertex to every neighbor of its corresponding original vertex, and adding a root vertex adjacent to all shadow vertices. Balakrishnan and Francis Raj \cite{balakrishnan} proved that, for every connected graph $G$ of order at least two, its Mycielskian graph $\mu(G)$ is maximally edge-connected. Subsequently, Guo, Liu and Guo \cite{guo} further proved that $\mu(G)$ is super edge-connected if and only if $G$ is not isomorphic to $K_2$.

In this paper, we focus on the restricted edge-connectivity of $\mu(G)$ and prove
that $\lambda'(\mu(G))=\min\{\xi(\mu(G)),\beta^{*}(G)\}$, where $\beta^{*}(G)$
is the minimum of $3|\partial_G(X)|+|X|$ over all restricted edge-cuts
$\partial_G(X)$ of $G$ and is set to $+\infty$ if $G$ contains no restricted edge-cut.

The remainder of this paper is organized as follows. Section~\ref{sec:prelim} introduces the notation and preliminary results. Section~\ref{sec:main} establishes the main theorem and its consequences. Section~\ref{sec:conclusion} concludes the paper.

\section{Preliminaries}\label{sec:prelim}

For terminology and notation not defined here, we follow Bondy and Murty
\cite{bondy}. All graphs considered in this paper are finite, simple and
undirected. Unless otherwise stated, $G=(V(G),E(G))$ denotes a connected graph
of order at least 2.

For a vertex $u\in V(G)$, let $N_G(u)$ denote the neighborhood of $u$ in $G$ and let
$d_G(u)=|N_G(u)|$ be its degree. The minimum degree
of $G$ is denoted by $\delta(G)$. For an edge $e=xy\in E(G)$, its edge-degree is
\(\xi_G(e)=d_G(x)+d_G(y)-2\). The minimum edge-degree of $G$ is defined as
 \(\xi(G)=\min_{e\in E(G)}\xi_G(e)\).

For a vertex subset $X\subseteq V(G)$, let $G[X]$ denote the subgraph induced by $X$ and $e(X)=|E(G[X])|$. For
$x\in V(G)$, define $d_X(x)=|N_G(x)\cap X|$. For disjoint vertex
subsets $X,Y\subseteq V(G)$, let $E(X,Y)=\{xy\in E(G):x\in X,\ y\in Y\}$ and
$e(X,Y)=|E(X,Y)|$. For $X\subseteq V(G)$, write $\overline {X}=V(G)\setminus X$ and
denote the edge boundary of $X$ by $\partial_G(X)=E(X, \overline {X})$.

We identify $V(G)$ with $[n]=\{1,\ldots,n\}$. For each $i\in[n]$, let $v_i$ and $u_i$ denote the original vertex and the shadow vertex of $\mu(G)$ corresponding to $i$, respectively, and let $w$ be the root. The Mycielskian graph $\mu(G)$ is defined as the graph with vertex set
\[
V(\mu(G))=\{v_i:i\in[n]\}\cup\{u_i:i\in[n]\}\cup\{w\}
\]
and edge set
\[
E(\mu(G))=\{v_iv_j:ij\in E(G)\}\cup\{v_iu_j,v_ju_i:ij\in E(G)\}\cup\{u_iw:i\in[n]\}.
\]
By the definition of the Mycielskian graph, $d_{\mu(G)}(v_i)=2d_G(i)$, $d_{\mu(G)}(u_i)=d_G(i)+1$ and $d_{\mu(G)}(w)=n$.

For $X\subseteq[n]$, let
\(V_X=\{v_i:i\in X\}\), \(U_X=\{u_i:i\in X\}\) and
\(S_X=V_X\cup U_X\). By the definition of $\mu(G)$, each edge
$ij\in\partial_G(X)$ with $i\in X$ and $j\in\overline X$ contributes exactly
the three edges $v_iv_j$, $v_iu_j$ and $v_ju_i$ to
$\partial_{\mu(G)}(S_X)$, while the root edges in
$\partial_{\mu(G)}(S_X)$ are precisely $u_iw$ with $i\in X$. Hence
\begin{equation}\label{eq:paired}
|\partial_{\mu(G)}(S_X)|
=
3|\partial_G(X)|+|X|.
\end{equation}

Let $\emptyset\ne X\subsetneq[n]$. By the definition of the Mycielskian graph, for each
$i\in X$, the vertices $v_i$ and $u_i$ have a neighbor in $S_X$ if and only if
$d_X(i)>0$. Hence $\mu(G)[S_X]$ has no isolated vertices if and only if
$G[X]$ has no isolated vertices. The complementary induced subgraph has vertex set
$V(\mu(G))\setminus S_X=V_{\overline X}\cup U_{\overline X}\cup\{w\}$.
Every vertex of $U_{\overline X}$ is adjacent to $w$ in the complementary induced
subgraph and $w$ is nonisolated because $\overline X\ne\emptyset$. By the definition of the Mycielskian graph, for each
$j\in\overline X$, the vertex $v_j$ has a neighbor in the complementary
induced subgraph if and only if $d_{\overline X}(j)>0$. Hence the
complementary induced subgraph has no isolated vertices if and only if
$G[\overline X]$ has no isolated vertices. 
If $G$ has a restricted edge-cut, define
\begin{equation}\label{eq:beta}
\beta^*(G)
=
\min
\left\{
3|\partial_G(X)|+|X|
:
\partial_G(X)\text{ is a restricted edge-cut of }G
\right\}.
\end{equation}
Otherwise, set $\beta^*(G)=+\infty$.

The following lemma determines the minimum edge-degree of $\mu(G)$.

\begin{lem}\label{lem:edge-degree}
Let $G$ be a connected graph of order $n\ge2$. Then
\[
\xi(\mu(G))
=\min\left\{\min_{ij\in E(G)}
\{2d_G(i)+d_G(j)-1,\,2d_G(j)+d_G(i)-1\},\,n+\delta(G)-1\right\}.
\]
\end{lem}
\begin{proof}
For $ij\in E(G)$, the degree formulas give
$\xi_{\mu(G)}(v_iv_j)=2d_G(i)+2d_G(j)-2$,
$\xi_{\mu(G)}(v_iu_j)=2d_G(i)+d_G(j)-1$ and
$\xi_{\mu(G)}(v_ju_i)=2d_G(j)+d_G(i)-1$.
Since $G$ is connected and $n\ge2$, we have $d_G(i),d_G(j)\ge1$,
which gives $\xi_{\mu(G)}(v_iv_j)\ge\xi_{\mu(G)}(v_iu_j)$ and
$\xi_{\mu(G)}(v_iv_j)\ge\xi_{\mu(G)}(v_ju_i)$.
For a root edge $u_kw$, we have $\xi_{\mu(G)}(u_kw)=n+d_G(k)-1$
and therefore $\min_{k\in[n]}\xi_{\mu(G)}(u_kw)=n+\delta(G)-1$.
Since every edge of $\mu(G)$ is one of these types,  the stated formula   follows.
\end{proof}

\begin{lem}\label{lem:upper}
Let $G$ be a connected graph of order $n\ge2$ and let $\mu(G)$ be its
Mycielskian graph. Then $\mu(G)$ has a restricted edge-cut and
\[
\lambda'(\mu(G))
\le
\min\{\xi(\mu(G)),\beta^*(G)\}.
\]
\end{lem}
\begin{proof}
Since $G$ is connected and $n\ge2$, the Mycielskian graph $\mu(G)$ is connected
and has order $2n+1\ge5$. Also, $d_G(i)\ge1$ for every $i\in[n]$,
which gives $\delta(\mu(G))\ge2$ by the degree formulas above.
Thus $\mu(G)$ is not a star and the result of Esfahanian and Hakimi
\cite{esfahanian} shows that it admits a restricted edge-cut and satisfies
$\lambda'(\mu(G))\le\xi(\mu(G))$.
If $\beta^*(G)=+\infty$, then the desired inequality follows from
$\lambda'(\mu(G))\le\xi(\mu(G))$. If $\beta^*(G)<+\infty$, let
$\partial_G(X)$ be a restricted edge-cut of $G$ satisfying
$3|\partial_G(X)|+|X|=\beta^*(G)$. By the definition of $\mu(G)$,
$\partial_{\mu(G)}(S_X)$ is a restricted edge-cut of size
$3|\partial_G(X)|+|X|=\beta^*(G)$, which gives
$\lambda'(\mu(G))\le\beta^*(G)$ and completes the proof.
\end{proof}

\section{Main theorem}\label{sec:main}

The following lemma describes a minimum restricted edge-cut when
$\mu(G)$ is not maximally restricted edge-connected.

\begin{lem}\label{lem:structural}
Let $G$ be a connected graph of order $n\ge2$ and let $\mu(G)$ be its
Mycielskian graph. If $\lambda'(\mu(G))<\xi(\mu(G))$, then there exists
$\emptyset\ne X\subsetneq V(G)$ such that
$\partial_{\mu(G)}(S_X)$ is a minimum restricted edge-cut of $\mu(G)$.
\end{lem}

\begin{proof}
Suppose, to the contrary, that there is no nonempty proper subset
$X\subseteq V(G)$ such that $\partial_{\mu(G)}(S_X)$ is a minimum
restricted edge-cut of $\mu(G)$.
Among all minimum restricted edge-cuts of $\mu(G)$, let $F$ be one for which
the component of $\mu(G)-F$ not containing $w$ has minimum order.
Let $H_1$ and $H_2$ be the two components of $\mu(G)-F$, where
$w\in V(H_2)$. And put $A=V(H_1)$ and $B=V(H_2)$.
Then $F=\partial_{\mu(G)}(A)$ and
$|F|=\lambda'(\mu(G))<\xi(\mu(G))\le\xi_{\mu(G)}(e)$ for every
$e\in E(\mu(G))$. By the definition of the Mycielskian graph and the
connectedness of $G$, $\delta(\mu(G))\ge2$.

If $|A|=2$, then by the connectedness of $\mu(G)[A]$, its two vertices
form an edge whose edge-degree equals $|F|$, contradicting
$\lambda'(\mu(G))<\xi(\mu(G))$. Thus $|A|\ge3$.

Suppose that a vertex $x\in A$ has only one neighbor in $A$. Since
$\mu(G)[A]$ is connected and $x$ has only one neighbor in $A$,
$\mu(G)[A\setminus\{x\}]$ is connected and has order at least two.
As $d_{\mu(G)}(x)\ge2$, $x$ has a neighbor in $B$ and
$\mu(G)[B\cup\{x\}]$ is connected. Hence $\partial_{\mu(G)}(A\setminus\{x\})$ is a restricted edge-cut and its size satisfies
$|\partial_{\mu(G)}(A\setminus\{x\})|=|F|+|N_{\mu(G)}(x)\cap A|-|N_{\mu(G)}(x)\cap B|=|F|+1-\bigl(d_{\mu(G)}(x)-1\bigr)=|F|+2-d_{\mu(G)}(x)\le |F|$.
Since it is a restricted edge-cut and $|F|=\lambda'(\mu(G))$, its size must equal $|F|$. Thus it is also a minimum restricted edge-cut. Its component not containing $w$ is $A\setminus\{x\}$, which has smaller order than $A$, contradicting the choice of $F$. Therefore
$\delta(\mu(G)[A])\ge2$.

Partition $V(G)$ according to the positions of the corresponding
original and shadow vertices:
\[
I=\{i:v_i,u_i\in A\},\enspace P=\{i:v_i\in A,u_i\in B\},\enspace
Q=\{i:u_i\in A,v_i\in B\},\enspace R=\{i:v_i,u_i\in B\}.
\]
Then $A=V_I\cup U_I\cup V_P\cup U_Q$ and
$B=U_P\cup V_Q\cup V_R\cup U_R\cup\{w\}$.

Each edge counted by $e(I,P)$, $e(I,Q)$, $e(I,R)$, $e(P)$, $e(P,Q)$,
$e(P,R)$, $e(Q)$ or $e(Q,R)$ contributes respectively $1,2,3,2,1,2,2$
or $1$ edges to $F$. Edges within $I$ or $R$ contribute none and
$u_iw\in F$ exactly when $i\in I\cup Q$. Hence
\[
\begin{aligned}
|F|
&=e(I,P)+2e(I,Q)+3e(I,R)+2e(P)+e(P,Q)\\
&\quad+2e(P,R)+2e(Q)+e(Q,R)+|I|+|Q|.
\end{aligned}
\]

We distinguish two cases according to whether $Q$ is empty or not empty.

\Needspace{5\baselineskip}
\smallskip
\noindent\textbf{Case 1. $Q\ne\emptyset$.}

We first prove that $e(P,Q)=0$. Otherwise, suppose that $e(P,Q)>0$. Let
$p_1q_1\in E(G)$, where $p_1\in P$ and $q_1\in Q$.
Then $v_{p_1}u_{q_1}\in E(\mu(G))$. Since $p_1q_1$ is counted in both $d_Q(p_1)$ and $d_P(q_1)$,
$d_Q(p_1)+d_P(q_1)-1\le e(P,Q)$. Combining this with
$d_I(p_1)\le|I|$, $d_Q(p_1)\le|Q|$ and the degree and edge counts, we obtain
\[
\begin{aligned}
\xi_{\mu(G)}(v_{p_1}u_{q_1})
&=2d_G(p_1)+d_G(q_1)-1\\
&=2\bigl(d_I(p_1)+d_P(p_1)+d_Q(p_1)+d_R(p_1)\bigr)\\
&\quad+d_I(q_1)+d_P(q_1)+d_Q(q_1)+d_R(q_1)-1\\
&=2d_I(p_1)+d_I(q_1)+2d_P(p_1)+d_Q(p_1)\\
&\quad+\bigl(d_Q(p_1)+d_P(q_1)-1\bigr)+d_Q(q_1)+2d_R(p_1)+d_R(q_1)\\
&\le e(I,P)+d_I(p_1)+e(I,Q)+2e(P)+d_Q(p_1)\\
&\quad+e(P,Q)+e(Q)+2e(P,R)+e(Q,R)\\
&\le e(I,P)+|I|+e(I,Q)+2e(P)+|Q|\\
&\quad+e(P,Q)+e(Q)+2e(P,R)+e(Q,R)\\
&\le e(I,P)+|I|+e(I,Q)+2e(P)+|Q|+e(P,Q)\\
&\quad+e(Q)+2e(P,R)+e(Q,R)+e(I,Q)+3e(I,R)+e(Q)\\
&=e(I,P)+2e(I,Q)+3e(I,R)+2e(P)\\
&\quad+e(P,Q)+2e(P,R)+2e(Q)+e(Q,R)+|I|+|Q|\\
&=|F|.
\end{aligned}
\]
a contradiction. Thus $e(P,Q)=0$.

We next prove that $|Q|=1$. Otherwise, suppose that $|Q|\ge2$. Let $q_2\in Q$ satisfy
$d_G(q_2)=\min\{d_G(q):q\in Q\}$. By the choice of $q_2$ and $e(P,Q)=0$,
\[
d_G(q_2)
\le\frac{e(I,Q)+2e(Q)+e(Q,R)}{|Q|}
\le\frac{e(I,Q)+2e(Q)+e(Q,R)}2.
\]
Since $\delta(\mu(G)[A])\ge2$, $u_{q_2}$ has a neighbor in $A$ and
$e(P,Q)=0$ forces every such neighbor to lie in $V_I$. Let
$i_1\in N_G(q_2)\cap I$. Then $v_{q_2}u_{i_1}\in E(\mu(G))$ and
$d_G(i_1)\le|I|-1+e(I,P)+e(I,Q)+e(I,R)$.
It follows that
\[
\begin{aligned}
\xi_{\mu(G)}(v_{q_2}u_{i_1})
&=2d_G(q_2)+d_G(i_1)-1\\
&\le \bigl(e(I,Q)+2e(Q)+e(Q,R)\bigr)\\
&\quad+\bigl(|I|-1+e(I,P)+e(I,Q)+e(I,R)\bigr)-1\\
&=e(I,P)+2e(I,Q)+e(I,R)+2e(Q)\\
&\quad+e(Q,R)+|I|-2\\
&\le e(I,P)+2e(I,Q)+e(I,R)+2e(Q)+e(Q,R)+|I|-2\\
&\quad+2e(I,R)+2e(P)+2e(P,R)+|Q|\\
&=e(I,P)+2e(I,Q)+3e(I,R)+2e(P)\\
&\quad+2e(P,R)+2e(Q)+e(Q,R)+|I|+|Q|-2\\
&=|F|-2\\
&<|F|.
\end{aligned}
\]
a contradiction. Therefore $|Q|=1$ and we write $Q=\{q\}$.

Having established $Q=\{q\}$, we prove that $d_I(q)\ge d_R(q)+2$. The inequality $\delta(\mu(G)[A])\ge2$ ensures that $\mu(G)[A\setminus\{u_q\}]$ has no isolated vertices. The subgraph $\mu(G)[B]$ is connected and $u_q$ is adjacent to $w\in B$. Hence $\mu(G)[B\cup\{u_q\}]$ is connected. Thus $\partial_{\mu(G)}(A\setminus\{u_q\})$ is a restricted edge-cut. By the minimality of $F$ and the choice of $A$, we have $|\partial_{\mu(G)}(A\setminus\{u_q\})|>|F|$. The conditions $Q=\{q\}$ and $e(P,Q)=0$ imply that every neighbor of $q$ in $G$ belongs to $I\cup R$. By the construction of the Mycielskian graph, $u_q$ has $d_I(q)$ neighbors in $A$ and $d_R(q)+1$ neighbors in $B$. Consequently, $|F|<|\partial_{\mu(G)}(A\setminus\{u_q\})|=|F|+|N_{\mu(G)}(u_q)\cap A|-|N_{\mu(G)}(u_q)\cap B|=|F|+d_I(q)-d_R(q)-1$. Since $d_I(q)$ and $d_R(q)$ are integers, we obtain $d_I(q)\ge d_R(q)+2$.

We then prove that $\mu(G)[B\setminus\{v_q\}]$ contains an isolated vertex. Suppose, to the contrary, that $\mu(G)[B\setminus\{v_q\}]$ has no isolated vertices. Since $Q=\{q\}$ and $e(P,Q)=0$, every neighbor of $q$ in $G$ belongs to $I\cup R$. By the construction of the Mycielskian graph, each neighbor of $q$ in $I$ gives two neighbors of $v_q$ in $A$, while each neighbor of $q$ in $R$ gives two neighbors of $v_q$ in $B$. Hence $|N_{\mu(G)}(v_q)\cap A|=2d_I(q)$ and $|N_{\mu(G)}(v_q)\cap B|=2d_R(q)$. Since $d_I(q)\ge d_R(q)+2$, we have $d_I(q)\ge2$. Thus $\mu(G)[A\cup\{v_q\}]$ has no isolated vertices. By our contrary assumption, $\mu(G)[B\setminus\{v_q\}]$ also has no isolated vertices. Hence $\partial_{\mu(G)}(A\cup\{v_q\})$ is a restricted edge-cut. Its size satisfies $|\partial_{\mu(G)}(A\cup\{v_q\})|=|F|+|N_{\mu(G)}(v_q)\cap B|-|N_{\mu(G)}(v_q)\cap A|=|F|+2d_R(q)-2d_I(q)\le|F|-4<|F|$, contradicting the minimality of $F$. Therefore, $\mu(G)[B\setminus\{v_q\}]$ contains an isolated vertex.

Every shadow vertex in $B\setminus\{v_q\}$ is adjacent to the root $w\in B\setminus\{v_q\}$. Hence the isolated vertex must be an original vertex $v_{r_1}\in B\setminus\{v_q\}$, where $r_1\in R$ by the definition of the partition. Since $\mu(G)[B]$ has no isolated vertices, the only neighbor of $v_{r_1}$ in $B$ is $v_q$. By the construction of the Mycielskian graph, this gives $d_P(r_1)=d_R(r_1)=0$ and $qr_1\in E(G)$. Consequently, $d_G(r_1)-1=d_I(r_1)\le e(I,R)$. Moreover, $v_qu_{r_1}\in E(\mu(G))$.

We finally derive a contradiction by calculating the edge-degree of $v_qu_{r_1}$. Since $Q=\{q\}$ and $e(P,Q)=0$, we have $d_G(q)=e(I,Q)+e(Q,R)$. Together with $d_G(r_1)-1=d_I(r_1)\le e(I,R)$, the calculation gives
\[
\begin{aligned}
\xi_{\mu(G)}(v_qu_{r_1})
&=2d_G(q)+d_G(r_1)-1\\
&=2(e(I,Q)+e(Q,R))+d_I(r_1)\\
&\le 2e(I,Q)+2e(Q,R)+e(I,R)\\
&<2e(I,Q)+2e(Q,R)+e(I,R)+e(I,P)+2e(I,R)+2e(P)\\
&\quad+2e(P,R)+|I|-e(I,Q)+e(I,Q)-e(Q,R)+1\\
&=e(I,P)+2e(I,Q)+3e(I,R)+2e(P)+2e(P,R)+e(Q,R)+|I|+1\\
&=|F|,
\end{aligned}
\]
where the strict inequality follows from $|I|-e(I,Q)\ge0$ and $e(I,Q)-e(Q,R)+1\ge3$. This is a contradiction.

\smallskip
\noindent\textbf{Case 2. $Q=\emptyset$.}

If $P=\emptyset$, then $A=S_I$ and $\partial_{\mu(G)}(S_I)$ is a
minimum restricted edge-cut, contrary to our assumption. Thus $P\ne\emptyset$
and $A=S_I\cup V_P$, $B=U_P\cup V_R\cup U_R\cup\{w\}$ and
$S_{I\cup P}=A\cup U_P$.

\Needspace{10\baselineskip}
We first prove that $R\ne\emptyset$. Otherwise, suppose that $R=\emptyset$. Then
$n=|I|+|P|$ and $|F|=e(I,P)+2e(P)+|I|=|I|+\sum_{p\in P}d_G(p)$.
Let $p_2\in P$ satisfy $d_G(p_2)=\min\{d_G(p):p\in P\}$. Since
$d_G(p_2)\ge1$, we have
$\xi_{\mu(G)}(u_{p_2}w)=n+d_G(p_2)-1\allowbreak
=|I|+|P|+d_G(p_2)-1\allowbreak
=|I|+|P|d_G(p_2)-(|P|-1)(d_G(p_2)-1)\allowbreak
\le |I|+|P|d_G(p_2)\allowbreak
=|I|+\sum_{p\in P}d_G(p_2)\allowbreak
\le |I|+\sum_{p\in P}d_G(p)=|F|$, contradicting
$|F|<\xi_{\mu(G)}(u_{p_2}w)$. Hence $R\ne\emptyset$.

We next prove that $d_I(p)\ge d_R(p)+1$ for every $p\in P$.
Otherwise, suppose that $d_I(p)\le d_R(p)$ for some $p\in P$.
Since $|N_{\mu(G)}(v_p)\cap B|=d_P(p)+2d_R(p)
\ge d_P(p)+2d_I(p)=|N_{\mu(G)}(v_p)\cap A|\ge2$,
$v_p$ has a neighbor in $B$ and $\mu(G)[B\cup\{v_p\}]$ is connected.
Since $\delta(\mu(G)[A])\ge2$, $\mu(G)[A\setminus\{v_p\}]$ has no
isolated vertices. Hence $\partial_{\mu(G)}(A\setminus\{v_p\})$ is a
restricted edge-cut and
$|\partial_{\mu(G)}(A\setminus\{v_p\})|
=|F|+\bigl(2d_I(p)+d_P(p)\bigr)-\bigl(d_P(p)+2d_R(p)\bigr)
=|F|+2\bigl(d_I(p)-d_R(p)\bigr)\le |F|$.
Since it is a restricted edge-cut and $|F|=\lambda'(\mu(G))$,
its size must equal $|F|$. Every component of
$\mu(G)[A\setminus\{v_p\}]$ has order smaller than $|A|$,
contradicting the choice of $F$. Therefore, by integrality,
$d_I(p)\ge d_R(p)+1$ for every $p\in P$.

We then prove that $\partial_{\mu(G)}(S_{I\cup P})$ is a restricted edge-cut.
Since $R\ne\emptyset$ and $I\cup P\ne\emptyset$, the paired set
$S_{I\cup P}$ is a nonempty proper subset of $V(\mu(G))$.
Since $A=V_{I\cup P}\cup U_I$ and $\delta(\mu(G)[A])\ge2$,
the definition of the Mycielskian graph gives, for every $i\in I\cup P$, a vertex
$j\in I\cup P$ with $ij\in E(G)$. Thus both $v_i$ and $u_i$ are
adjacent to $v_j\in S_{I\cup P}$ and $\mu(G)[S_{I\cup P}]$ has no
isolated vertices.

\Needspace{6\baselineskip}
The induced subgraph $\mu(G)[V(\mu(G))\setminus S_{I\cup P}]$ has vertex set
$V_R\cup U_R\cup\{w\}$. Since $R\ne\emptyset$ and every vertex of $U_R$
is adjacent to $w$, only vertices of $V_R$ can be isolated. Suppose that
$v_{r_2}$, $r_2\in R$, is isolated. Then $d_R(r_2)=0$. As $\mu(G)[B]$
has no isolated vertices, $v_{r_2}$ has a neighbor in $U_P$.
Let $p_3\in N_G(r_2)\cap P$ have minimum degree among the vertices of
$N_G(r_2)\cap P$. \mbox{Then $v_{r_2}u_{p_3}\in E(\mu(G))$ and $d_G(p_3)\le\frac{\sum_{p\in N_G(r_2)\cap P}d_G(p)}{d_P(r_2)}\le\frac{\sum_{p\in P}d_G(p)}{d_P(r_2)}$}.
Since $\sum_{p\in P}d_G(p)\ge e(P,R)\ge d_P(r_2)\ge1$, we obtain
\[
\begin{aligned}
\xi_{\mu(G)}(v_{r_2}u_{p_3})
&=2d_G(r_2)+d_G(p_3)-1\\
&=2d_I(r_2)+2d_P(r_2)+d_G(p_3)-1\\
&\le2e(I,R)+2d_P(r_2)+\frac{\sum_{p\in P}d_G(p)}{d_P(r_2)}-1\\
&=2e(I,R)+2d_P(r_2)+\frac{\sum_{p\in P}d_G(p)-d_P(r_2)}{d_P(r_2)}\\
&\le2e(I,R)+\sum_{p\in P}d_G(p)+d_P(r_2)\\
&=2e(I,R)+e(I,P)+2e(P)+e(P,R)+d_P(r_2)\\
&\le e(I,P)+2e(I,R)+2e(P)+2e(P,R)\\
&\le e(I,P)+3e(I,R)+2e(P)+2e(P,R)+|I|\\
&=|F|,
\end{aligned}
\]
a contradiction. Therefore the complementary induced subgraph has no
isolated vertices and $\partial_{\mu(G)}(S_{I\cup P})$ is a restricted edge-cut.

Finally, we will prove that $|\partial_{\mu(G)}(S_{I\cup P})|\le|F|$.
By the definition of the Mycielskian graph,
\[
|\partial_{\mu(G)}(S_{I\cup P})|
=3\bigl(e(I,R)+e(P,R)\bigr)+|I|+|P|.
\]
Since $d_I(p)\ge d_R(p)+1$ for every $p\in P$, we have
$e(I,P)=\sum_{p\in P}d_I(p)\ge\sum_{p\in P}d_R(p)+|P|=e(P,R)+|P|$.
Combining these two estimates, we obtain
\[
\begin{aligned}
|\partial_{\mu(G)}(S_{I\cup P})|
&=3e(I,R)+3e(P,R)+|I|+|P|\\
&=3e(I,R)+2e(P,R)+|I|+\bigl(e(P,R)+|P|\bigr)\\
&\le3e(I,R)+2e(P,R)+|I|+e(I,P)\\
&\le3e(I,R)+2e(P,R)+|I|+e(I,P)+2e(P)\\
&=e(I,P)+3e(I,R)+2e(P)+2e(P,R)+|I|\\
&=|F|,
\end{aligned}
\]
Since $\partial_{\mu(G)}(S_{I\cup P})$ is a restricted edge-cut and
$|F|=\lambda'(\mu(G))$, its size is at least $|F|$.
Thus equality holds, so it is a minimum restricted edge-cut of $\mu(G)$,
contrary to our assumption.
\end{proof}

%
\begin{thm}\label{thm:main}
Let $G$ be a connected graph of order $n\ge2$. Then
\[
\lambda'(\mu(G))=\min\{\xi(\mu(G)),\beta^*(G)\}.
\]
\end{thm}
\begin{proof}
Lemma~\ref{lem:upper} gives
$\lambda'(\mu(G))\le\min\{\xi(\mu(G)),\beta^*(G)\}$.
It remains to prove $\lambda'(\mu(G))\geq\min\{\xi(\mu(G)),\beta^*(G)\}$.

If $\lambda'(\mu(G))=\xi(\mu(G))$, then $\lambda'(\mu(G))\geq\min\{\xi(\mu(G)),\beta^*(G)\}$ holds.
Suppose that $\lambda'(\mu(G))<\xi(\mu(G))$.
By Lemma~\ref{lem:structural}, there exists a nonempty proper subset
$X\subseteq V(G)$ such that $\partial_{\mu(G)}(S_X)$ is a minimum
restricted edge-cut of $\mu(G)$.
The analysis of this restricted edge-cut in Section~\ref{sec:prelim} shows that
neither $G[X]$ nor $G[\overline X]$ has an isolated vertex.
Thus $\partial_G(X)$ is a restricted edge-cut of $G$.
By the definition of $\beta^*(G)$,
$\beta^*(G)\le3|\partial_G(X)|+|X|
=|\partial_{\mu(G)}(S_X)|=\lambda'(\mu(G))$.
Since $\lambda'(\mu(G))<\xi(\mu(G))$, we obtain
$\beta^*(G)<\xi(\mu(G))$.
Hence $\lambda'(\mu(G))\geq\beta^*(G)=\min\{\xi(\mu(G)),\beta^*(G)\}$.
\end{proof}

\begin{cor}\label{cor:iff}
Let $G$ be a connected graph of order $n\ge2$. Then $\mu(G)$ is maximally restricted edge-connected if and only if $\xi(\mu(G))\le\beta^*(G)$.
\end{cor}
\begin{proof}
By Theorem~\ref{thm:main}, $\lambda'(\mu(G))=\xi(\mu(G))$ if and only if
$\min\{\xi(\mu(G)),\beta^*(G)\}=\xi(\mu(G))$, which is equivalent to
$\xi(\mu(G))\le\beta^*(G)$.
\end{proof}

\begin{pro}\label{prop:bounds}
Let $G$ be a connected graph of order $n\ge2$. If $G$ has a restricted edge-cut, then
\(3\lambda'(G)+2\le\beta^*(G) \le3\lambda'(G)+\left\lfloor\frac n2\right\rfloor.\)
\end{pro}
\begin{proof}
For the lower bound, let $\partial_G(X)$ be any restricted edge-cut of $G$.
Since $|\partial_G(X)|\ge\lambda'(G)$ and $|X|\ge2$, we have
$3|\partial_G(X)|+|X|\ge3\lambda'(G)+2$.
Taking the minimum over all such $X$ gives
$\beta^*(G)\ge3\lambda'(G)+2$.

For the upper bound, let $F$ be a minimum restricted edge-cut of $G$ and
let $X$ be the vertex set of the smaller component of $G-F$.
Then $\lambda'(G)=|\partial_G(X)|$.
Since $|X|\le\lfloor n/2\rfloor$, it follows
$\beta^*(G)\le3|\partial_G(X)|+|X|\le3\lambda'(G)+\left\lfloor\frac n2\right\rfloor$.
\end{proof}

The following bounds on $\beta^*(G)$ yield sufficient conditions for $\mu(G)$ to be maximally restricted edge-connected.

\begin{cor}\label{cor:sufficient}
Let $G$ be a connected graph of order $n\ge2$.
If $G$ has a restricted edge-cut and $\xi(\mu(G))\le3\lambda'(G)+2$, then $\mu(G)$ is maximally restricted edge-connected. In particular, $\mu(G)$ is maximally restricted edge-connected whenever $n+\delta(G)-1\le3\lambda'(G)+2$.
\end{cor}
\begin{proof}
Since $\xi(\mu(G))\le3\lambda'(G)+2$, we obtain  
$\xi(\mu(G))\le3\lambda'(G)+2\le\beta^*(G)$ by Proposition~\ref{prop:bounds}, and   $\mu(G)$ is maximally restricted edge-connected by Corollary~\ref{cor:iff}.

If $n+\delta(G)-1\le3\lambda'(G)+2$, then,   by Lemma~\ref{lem:edge-degree},
$\xi(\mu(G))\le n+\delta(G)-1\le3\lambda'(G)+2$. Thus Corollary~\ref{cor:iff} implies $\mu(G)$ is maximally restricted edge-connected.
\end{proof}

\begin{cor}\label{cor:preservation}
Let $G$ be a  connected graph of order at least two. If $G$ has a restricted edge-cut and $\lambda'(G)\geq\frac{2}{3}\xi(G)$,
then $\mu(G)$ is maximally restricted edge-connected.
\end{cor}
\begin{proof}
By Corollary~\ref{cor:sufficient}, it suffices to prove
$\xi(\mu(G))\le3\lambda'(G)+2$.
Let $ij\in E(G)$ satisfy $d_G(i)+d_G(j)-2=\xi(G)$.
As $d_G(j)\ge1$, Lemma~\ref{lem:edge-degree} gives
\[
\begin{aligned}
\xi(\mu(G))
\le2d_G(i)+d_G(j)-1
\le2\bigl(d_G(i)+d_G(j)-2\bigr)+2=2\xi(G)+2\le3\lambda'(G)+2.
\end{aligned}
\]
\end{proof}

\begin{cor}\label{cor:preservation}
Let $G$ be a maximally restricted edge-connected graph of order at least two, that is, $\lambda'(G)=\xi(G)$.
Then $\mu(G)$ is maximally restricted edge-connected.
\end{cor}

\Needspace{15\baselineskip}
\section{Conclusions}\label{sec:conclusion}

In this paper, we determine the restricted edge-connectivity of the
Mycielskian graph $\mu(G)$ for every connected graph $G$ of order at least two
and establish the exact formula
$\lambda'(\mu(G))=\min\{\xi(\mu(G)),\beta^*(G)\}$.
We also show that, whenever $\lambda'(\mu(G))<\xi(\mu(G))$, a minimum
restricted edge-cut can be chosen as the boundary of a set consisting of
original vertices and their corresponding shadow vertices. Furthermore, we
establish several sufficient conditions for $\mu(G)$ to be maximally restricted edge-connected. Future work may
focus on $k$-restricted edge-connectivity of Mycielskian graphs for $k\geq 3$
and on extending the present result to generalized Mycielskian graphs.


\end{document}